\documentclass[12pt, reqno]{amsart}
\usepackage{amssymb,amsthm,amsmath}
\usepackage[numbers,sort&compress]{natbib}
\usepackage{amssymb,amsmath}
\usepackage{amsfonts}
\usepackage{mathrsfs}
\usepackage{latexsym}
\usepackage{amssymb}
\usepackage{amsthm}
\usepackage{pdfsync}
\usepackage{indentfirst}
\date{Dec 10, 2019}

\newcommand{\R}{{\mathbb R}}

\renewcommand{\div}{{\rm div}}

\let\oldsection\section
\renewcommand\section{\setcounter{equation}{0}\oldsection}

\newtheorem{theorem}{Theorem}[section]
\newtheorem{lemma}{Lemma}[section]
\newtheorem{proposition}{Proposition}[section]

\newtheorem{remark}{Remark}[section]

\def\ba{\begin{eqnarray}}
\def\ea{\end{eqnarray}}

\def\R{\Bbb R}

\newcommand{\beq}{\begin{equation}}
\newcommand{\eeq}{\end{equation}}
\newcommand{\ben}{\begin{eqnarray}}
\newcommand{\een}{\end{eqnarray}}
\newcommand{\beno}{\begin{eqnarray*}}
\newcommand{\eeno}{\end{eqnarray*}}

\allowdisplaybreaks
\begin{document}

\title[Asymptotic properties of the plane shear thickening fluids]{Asymptotic properties of the plane shear thickening fluids with bounded energy integral}

\author{Shuai~Li}
\address[Shuai~Li]{School of Mathematical Sciences, Dalian University of Technology, Dalian, 116024,  China}
\email{ls21701048@163.com}

\author{Tao~Wang}
\address[Tao~Wang]{School of Mathematical Sciences, Dalian University of Technology, Dalian, 116024,  China}
\email{muyuwt@163.com}

\author{Wendong~Wang}
\address[Wendong~Wang]{School of Mathematical Sciences, Dalian University of Technology, Dalian, 116024,  China}
\email{wendong@dlut.edu.cn}




\begin{abstract}
In this note we investigate the asymptotic behavior of plane shear thickening fluids around a bounded obstacle. Different from the Navier-Stokes case considered by Gilbarg-Weinberger in \cite{GW1978}, where the good structure of the vorticity can be exploited and weighted energy estimates can be applied, we have to overcome the nonlinear term of high order. The decay estimates of the velocity was obtained by combining Point-wise Behavior Theorem in \cite{Galdi} and Brezis-Gallouet inequality in \cite{BG1980} together, which is independent of interest.

\end{abstract}

\maketitle

\allowdisplaybreaks

\section{Introduction}\label{sec1}

As Ladyzhenskaya suggests in her monograph in \cite{Lady1969}, it is interesting  to investigate "new equations for the
description of the motion of viscous incompressible fluids", which roughly speaking means to consider
viscosity coefficients, which depend on the modulus of the symmetric gradient,
\beno
\varepsilon(u) = \frac{1}{2} (Du + (Du)^T) = \frac{1}{2} (\partial_iu_k + \partial_ku_i)_{1 \leq i,k \leq 2}
\eeno
of the velocity field $u$, for example, in a monotonically increasing way (shear thickening
case). In this note we will consider this problem in a very special situation restricting ourselves to
stationary flows through an exterior domain $\Omega \subset \mathbb{R}^2$ with smooth boundary $\partial\Omega$. More precisely,
consider the solution $u:\Omega\rightarrow\mathbb{R}^2,\pi:\Omega\rightarrow\mathbb{R}$ of the following system
\begin{eqnarray}\label{eq:gns}
 \left\{
    \begin{array}{llll}
    \displaystyle -\div[T(\varepsilon(u))]+u_k\partial_ku+D\pi=0,~{\rm in}~\Omega,\\
    \displaystyle \div~u=0,~{\rm in}~\Omega,\\
    \end{array}
 \right.
\end{eqnarray}
where the $\Omega = {\R}^2 \setminus \overline{B_{R_0}(0)}$. More details on viscous incompressible flow, we refer to \cite{Lady1969,Galdi1994-1,Galdi1994-2,FS2000,FMS2003,Beirao2005, Fu2012exist} and the references therein.

The system (\ref{eq:gns}) describes the stationary flow of an incompressible generalized Newtonian fluid, where $u$ is the velocity field, $\pi$ is the pressure function, $u_k\partial_ku$ is the convective term, and $T$ represents the stress deviator tensor.
 And we use $\varepsilon(u)$ to stand for the symmetric part of the differential matrix $Du$ of $u$.
We assume that the stress tensor $T$ is the gradient of a potential $H : S^{2 \times 2} \rightarrow \mathbb{R}$ defined on the space $S^{2 \times 2}$ of all symmetric $2 \times 2$ matrices of the following form $$H(\varepsilon)= h(|\varepsilon|),$$ where $h$ is a nonnegative function of class $C^3$.
Thus
\ben\label{eq:h defination}
T(\varepsilon) = DH(\varepsilon) = \mu(|\varepsilon|)\varepsilon , \  \mu(t) = \frac{h^{'}(t)}{t}.
\een
Note that the Navier-Stokes equations for incompressible Newtonian fluids  follow from the system (\ref{eq:gns})  if  $\mu$ is a constant.
If $\mu$ is not a constant, it means that the viscosity coefficient depends on $\varepsilon$, and system (\ref{eq:gns}) describes the motion of continuous media of generalized Newtonian fluids.

As in \cite{Fu2012Liou}, assume that the potential $h$ satisfies the follow conditions:
\begin{center}
$h$ is strictly increasing and convex \phantom{hfill}\\
\hfill together with $h^{''}(0)> 0$ and $\lim\limits_{t \rightarrow 0} \frac{h(t)}{t} = 0$ ; \hfill (A1)\\
\hfill there exist a constant $a \ge 1$ such that $h(2t) \leq ah(t)$ for all $t \geq 0$; \hfill (A2)\\
\hfill $\frac{h^{'}(t)}{t} \leq h^{''}(t)$ for any $t \geq 0$. \hfill (A3)\\
\end{center}

Let us sketch some progress on the system (\ref{eq:gns}). First, the existence of strong solutions is proved in a bounded domain by Fuchs in \cite{Fu2012exist}.
The existence of Dirichlet energy solutions satisfying the boundary condition at infinity in an exterior is very difficult, even if for the Navier-Stokes equations, which is related to Leray's question; for example see Leray \cite{Leray},  Amick \cite{Amick}, Russo \cite{Ru2009}, Pileckas-Russo \cite{PR2012} and the references therein. Different from the Navier-Stokes case, the regularity is also unknown for general form $h(t)$ as in \eqref{eq:h defination}. Bildhauer-Fuchs-Zhong \cite{BFZ2005} proved the weak solution is $C^{1,\alpha}$ by assuming $h(t)=t^2(1+t)^m$, see also recent improved result by Jin-Kang in \cite{JK2014}. The Liouville property of (\ref{eq:gns}) was
started by Fuchs in \cite{Fu2012Liou}, and later studied by Zhang in \cite{ZG2013} and \cite{ZG2015}, where they obtained the trivial property of the solution with the help of $u\in L^\infty$ or $\int_{\Omega}h(|\nabla u|)<\infty.$ The degenerate case $h(t)=t^p$ was also considered by Bildhauer-Fuchs-Zhang in \cite{BFZ2013} by assuming that $\int_{\Omega}|\nabla u|^p<\infty.$ More developments, we refer to \cite{Fu2012exist,BFZ2013} and the references therein.

In this note, motivated by the work of Gilbarg-Weinberger \cite{GW1978}, we investigate the asymptotic properties of the solutions of (\ref{eq:gns}).
In \cite{GW1978}, they showed that pressure function $\pi$ has a limit at infinity, $u(z) = o(\ln^{\frac12} r)$, and $|D u|\leq o(r^{-\frac34}(\ln r)^{\frac98})$ provided that the Dirichlet energy is bounded in an exterior domain, i.e., $\int_{\Omega} |Du|^2 dx < \infty$. Their proof relies on the fact that the vorticity of the 2D Navier-Stokes equations satisfies a nice elliptic equation,
to which the maximum principle applies. Here we consider
the case of shear thickening fluids, for $h$ satisfying the (A1)-(A3); however, it's difficult to exploit the good structure of the vorticity and apply weighted energy estimates, since the main part in (\ref{eq:gns}) is nonlinear. Inspired by  Point-wise Behavior Theorem in \cite{Galdi} and Brezis-Gallouet inequality (for example, see \cite{BG1980} or \cite{CPZ}, we obtain the higher energy estimates, which imply the decay estimates by combining point-wise behavior theorem and Brezis-Gallouet inequality together.

As show in \cite{BFZ2013}, the following properties of functions $h$ follows from (A1)-(A3).

(i) $\mu(t) = \frac{h^{'}(t)}{t}$ is an increasing function.

(ii)We have $h(0) = h^{'}(0) = 0$ and
\begin{eqnarray}\label{eq:h''>0}
    \begin{array}{llll}
    h(t) \geq \frac{1}{2} h^{''}(0)t^2,\quad h^{''}(0) > 0
    \end{array}
\end{eqnarray}

(iii) There exists a constant $a>0$ such that the function $h$ satisfies the balancing condition,
\begin{eqnarray}\label{eq:th'<h}
    \begin{array}{llll}
    \frac{1}{a}th^{'}(t) \leq h(t) \leq th^{'}(t),~~~~~t \ge 0
    \end{array}
\end{eqnarray}

(iv)
From the assumptions on $h$, we know the system satisfies the following elliptic condition, $\forall \varepsilon , \sigma \in S^2$,
\ben\label{eq:D2H LOWER BOUND}
    \frac{h^{'}(|\varepsilon|)}{|\varepsilon|} |\sigma|^2 \leq D^2H(\varepsilon)(\sigma,\sigma) \leq h^{''}(|\varepsilon|)|\sigma|^2,
\een
from which, together with (\ref{eq:h''>0}) and (\ref{eq:th'<h}), it follows that
\begin{eqnarray}\label{eq:D2H LOWER BOUND2}
    \begin{array}{llll}
    D^2H(\varepsilon)(\sigma,\sigma) \geq \frac12 h^{''}(0)|\sigma|^2.
    \end{array}
\end{eqnarray}

Let $T_r=B_r(0) \setminus \overline{B_{R_0}(0)}$ for any $r>R_0>0.$
Our first result is to estimate the $L^2$ norm of $D^2 u.$

\begin{proposition}\label{prop:nabla2 u}
 For $\Omega = \mathbb{R}^2 \setminus \overline{B_{R_0}(0)}$, let $u \in C^2(\Omega , \mathbb{R}^2) $ be a solutions of (\ref{eq:gns}).
 Then there hold $\int_\Omega | D^2u |^2  dx<\infty,$
\ben\label{eq:nabla2u decay}
\int_{T_{\frac{3}{2}r} \backslash T_{\frac{9}{8}r}}| D^2u |^2  dx \leq C\frac{\sqrt{\mathrm{log}r}}{r}, \een
and
\ben\label{eq:u decay}
|u(x)|= o(\sqrt{\ln(|x|)}),
\een
provided that $r=|x|$ is large enough and $\int_{\Omega} h(|D u|) dx <\infty$.
\end{proposition}

\begin{remark}
The $L^2$ norm of $D^2u$ was obtained in \cite{ZG2013} in the whole space $\mathbb{R}^2$,
\beno
\int_{B_R }| D^2u |^2  dx \leq C+C\frac{1}{R^3}\int_{B_{2R}}|u|^2dx
\eeno
Here we refine the estimate, especially for the exterior domain. One observation is to apply the Wirtinger's inequality in $L^3$ norm and another  technique is to use Poincar\'{e}-Sobolev inequality  in a circular region. The estimate \eqref{eq:u decay}
 is the same as the Navier-Stokes case in \cite{GW1978}.
 \end{remark}


In order to obtain the decay estimate of $Du$, we need to estimate the norm of $D^3u.$
Generally, it's more difficult. Motivated by the anisotropic variational problems in \cite{BF2003,ABF2005}, we give the following assumption.

(v) Assume that
\ben\label{eq;h'' lower bound}
1+C_1t^{\gamma_1} \leq h^{''}(t) \leq 1+C_2t^{\gamma_2},t \geq 0,
\een
where $\gamma_2\geq \gamma_1\geq 1$ and $C_2>C_1>0.$
From this, it follows that
\ben\label{eq:h' lower bound}
   &&t^2+C_1't^{2+\gamma_1} \leq  h(t) \leq t^2+C_2't^{2+\gamma_2},\nonumber\\
   &&t+C_1''t^{1+\gamma_1} \leq h^{'}(t) \leq t+C_2''t^{1+\gamma_2},\quad t \geq 0,
\een
where $C_1',C_2',C_1'',C_2''$ depend on $C_1,C_2,\gamma_1$ and $\gamma_2.$

(vi) Assume that
\ben\label{eq;h''' up bound}
h^{'''}(t) \leq C_3 + C_3't^{\gamma_0},\quad \gamma_0 \geq 0,\quad t \geq 0.
\een
where $C_3,C_3' \geq 0.$

\begin{proposition}\label{prop:nabla3 u}\ Let $\Omega = \mathbb{R}^2 \setminus \overline{B_{R_0}(0)}$,\ $u \in C^3(\Omega , \mathbb{R}^2)$ is  a solutions of (\ref{eq:gns}), if $|\varepsilon(u)| \leq C $, then there exists a constant $r_0>0$ such that
$$\int_{T_{\frac{3}{2}r} \setminus {T_{\frac{9}{8}r}}} | D^3u |^2  dx \leq  C ,$$
for any $r> r_0$, provided that $\int_\Omega h(|D u|)dx < \infty$.
\end{proposition}

\begin{remark} The conditions (v) and (vi) are used to estimate the term $J_3$ in Sec. 4:
\beno
 \int_{\Omega} \partial_n \left(\frac{h^{''}(|\varepsilon(u)|)|\varepsilon(u)| - h^{'}(|\varepsilon(u)|)}{|\varepsilon(u)|^3}\right) (\varepsilon(u) : \varepsilon(\partial_ku)) (\varepsilon(u) : \varepsilon(\partial_{kn}u)) \zeta^3 dx,
\eeno
which need the smallness of the error of $th''(t)-h'(t)$ and the growth control of $h'''(t)$.

 \end{remark}

With the help of  Proposition \ref{prop:nabla2 u} and \ref{prop:nabla3 u}, an argument by Brezis-Gallouet inequality yields that
\begin{theorem}\label{thm:main} Suppose that $u \in C^3(\Omega,\mathbb{R}^2)$ is a solution of (\ref{eq:gns}) and $ |\varepsilon(u)| \leq C,$ then there exists a constant $r_1>0$ such that \beno
||Du||_{L^{\infty}(T_{2r} \setminus {T_r})} \leq Cr^{-\frac{1}{2}}(\mathrm{log}r)^{\frac{3}{4}},
 \eeno
for all $r>r_1$, provided that $\int_{\Omega} h(|Du|) dx < \infty. $
\end{theorem}

Our paper is organized as follows: in Sec.2 we prove the proof of Theorem \ref{thm:main} under the assumptions of Proposition \ref{prop:nabla2 u} and \ref{prop:nabla3 u} by Brezis-Gallouet inequality.  In Sect.3 we are aimed to proving Proposition \ref{prop:nabla2 u} by Point-wise Behavior Theorem in \cite{Galdi} and in the last section we complete the proof of Proposition \ref{prop:nabla3 u}.

Throughout this paper we adopt the Einstein summation convention, which
means that the sum is taken with respect to indices repeated twice. Moreover,
throughout the remaining section, we denote by $C$ a general positive constant which depends only on known constant coefficients or norms and may be different from line to line.

\section{Proof of Theorem \ref{thm:main}}

\setcounter{equation}{0} \setcounter{theorem}{0}

Under the assumptions of Proposition \ref{prop:nabla2 u} and \ref{prop:nabla3 u}, we complete the proof of the main theorem.

\begin{proof}

Assume that $R_{0}=1$ without loss of generality, and define $$\widetilde{u}(\widetilde{x})=ru(r\widetilde{x})=ru(x),$$ where $x\in T_{2r} \setminus {T_{r}}$ and $r$ is large enough.
By Lemma \ref{lem:BG}, we have
$$\|D\widetilde{u}\|_{L^{\infty}(T_{2} \setminus  {T_{1}})}\leq C(1+\|D\widetilde{u}\|_{H^{1}(T_{2} \setminus  {T_{1}})}) \sqrt{{\mathrm{log}}(e+\|D^{3}\widetilde{u}\|_{L^{2}(T_{2} \setminus  {T_{1}})})},$$
And due to scaling, we have
$$\|D\widetilde{u}\|_{L^{\infty}(T_{2} \setminus  {T_{1}})}=r^{2}\|Du\|_{L^{\infty}(T_{2r} \setminus  {T_{r}})},$$
$$\|D\widetilde{u}\|_{L^{2}(T_{2} \setminus {T_{1}})}=r\|Du\|_{L^{2}(T_{2r} \setminus  {T_{r}})},$$
$$\|D^{2}\widetilde{u}\|_{L^{2}(T_{2} \setminus  {T_{1}})}=r^{2}\|D^{2}u\|_{L^{2}(T_{2r} \setminus  {T_{r}})},$$
$$\|D^{3}\widetilde{u}\|_{L^{2}(T_{2} \setminus {T_{1}})}=r^{3}\|D^{3}u\|_{L^{2}(T_{2r} \setminus  {T_{r}})}.$$
Hence
\ben\label{eq:B-G}
r^{2}\|Du\|_{L^{\infty}(T_{2r} \setminus  {T_{r}})}&\leq& C(1+r\|Du\|_{L^{2}(T_{2r} \setminus  {T_{r}})}+r^{2}\|D^{2}u\|_{L^{2}(T_{2r} \setminus {T_{r}})}) \nonumber\\
&&\cdot \sqrt{\mathrm{log}(e+r^{3}\|D^{3}u\|_{L^{2}(T_{2r} \setminus  {T_{r}})})},
\een
Due to the Proposition \ref{prop:nabla2 u}, there holds
$$\int_{T_{2r} \setminus  {T_{r}}}|D^{2}u|^{2}dx\leq \frac{C\sqrt{\mathrm{log}r}}{r}.$$
Besides, due to Proposition \ref{prop:nabla3 u}, we also get
$$\int_{T_{2r} \setminus  {T_{r}}}|D^{3}u|^{2}dx\leq C .$$
Then using (\ref{eq:B-G}), one can get
\begin{equation}
\begin{aligned}
r^{2}\|Du\|_{L^{\infty}(T_{2r} \setminus  {T_{r}})}& \leq C(1+r\|Du\|_{L^{2}(T_{2r} \setminus  {T_{r}})}+r^{2}\|D^{2}u\|_{L^{2}(T_{2r} \setminus {T_{r}})}) \\
&\mathrel{\phantom{=======}} \times \sqrt{\mathrm{log}(e+r^{3}\|D^{3}u\|_{L^{2}(T_{2r} \setminus  {T_{r}})})} \\
& \leq C(1+Cr+Cr^{\frac{3}{2}}(\mathrm{log}r)^{\frac{1}{4}})\sqrt{\mathrm{log}(e+Cr^{3})}\\
& \leq Cr^{\frac{3}{2}}(\mathrm{log}r)^{\frac{1}{4}}\sqrt{\mathrm{log}r},
\end{aligned}
\end{equation}
which implies
$$\|Du\|_{L^{\infty}(T_{2r} \setminus  {T_{r}})}\leq Cr^{-\frac{1}{2}}(\mathrm{log}r)^{\frac{3}{4}}.$$
The proof of Theorem \ref{thm:main} is complete.

\end{proof}

\section{Proof of Proposition \ref{prop:nabla2 u}}

\setcounter{equation}{0} \setcounter{theorem}{0}

Assume that $R_{0}=1$ without loss of generality and $\Omega=\mathbb{R}^2\setminus{\overline{B_1}}$. Moreover, $T_r=B_r\setminus{\overline{B_1}}$ for any $r>1.$

To prove the decay of $D^2u$, on one hand we explore the weak maximum principle of the equation (\ref{eq:gns}), which is similar to the vorticity form of the Navier-Stokes equations; on the other hand, an obstacle is to deal with the exterior domain and it seems that one can't apply the embedding theorem with scaling domain and the iterative lemma \ref{lem:Gia} simultaneously.
Our main idea is to apply the embedding theorem and  the iterative lemma \ref{lem:Gia} in different domains. More precisely, we use the embedding theorem in a circular domain with a proportional boundary via  weak maximum principle of the equation (\ref{eq:gns}), but  the iterative lemma is applies in the whole domain.

The proof is divided into two steps, at first we want to prove the $L^2$ norm of $D^2u$ is bounded and at last we prove the decay of this norm.

First, we introduce
the cut-off functions, which will be used in the next proof.

{\bf Case I. The choosing of $\phi$.} Assume that $\phi(x) \in C_0^{\infty}(\Omega)$ with $ 0 \le \phi \le 1$, satisfying that\\
i) for $r>10,$ $\rho>0$ and $\tau>0$, there holds $\frac{3}{4}r \le \rho < \tau \le r$;\\
ii)
\begin{equation}\label{eq:phi}\phi(x)=\phi(|x|)=
\left\{\begin{array}{llll}
1,\quad ~{\rm in}~T_{\rho} \backslash T_{3};\\
0, \quad ~{\rm in} ~T_{2},~{\rm and}~ \Omega\backslash T_{\tau}
\end{array}\right.
\end{equation}
iii)
\beno
&&|D\phi| \leq \frac{C}{\tau - \rho},\quad |D^2\phi| \leq \frac{C}{(\tau - \rho)^2},\quad {\rm as}\,\,x\in~ T_\tau\setminus{T_\rho};\\
&& |D\phi| \leq C,\quad |D^2\phi| \leq C,\quad {\rm as}\,\,x\in~ T_3\setminus{T_2}
\eeno

{\bf Case 2. The choosing of $\psi$.} Assume that $\psi(x) \in C_0^{\infty}(\Omega)$ with $ 0 \le \psi \le 1$, satisfying that\\
i) for $r\gg10,$ there holds
\begin{equation}\label{eq:psi}\psi(x)=\psi(|x|)=
\left\{\begin{array}{llll}
1,\quad ~{\rm in}~T_{2r} \backslash T_{r};\\
0, \quad ~{\rm in} ~T_{\frac{r}{2}},~{\rm and}~ \Omega\backslash T_{3r}
\end{array}\right.
\end{equation}
ii)
\beno
|D\psi| \leq \frac{C}{r},\quad |D^2\psi| \leq \frac{C}{r^2}.
\eeno

\begin{proof} First, we want to obtain Caccioppoli-type inequality by following the same route as in \cite{Fu2012Liou} or \cite{ZG2013}. Then we estimate the crucial items in more delicate analysis.

Choose the test function $\varphi_k = \partial_ku\,\eta^{2}$, where the cut-off function $\eta \in C_0^{\infty}(\Omega)$ with $ 0 \le \eta \le 1$. Multiply (\ref{eq:gns}) with $\partial_k\varphi_k$,  and integration by parts yields
\begin{equation}
\nonumber
\begin{aligned}
\int_{\Omega} \partial_k\sigma : \varepsilon(\varphi_k) dx - \int_{\Omega} D\pi \cdot \partial_k\varphi_k dx- \int_{\Omega} u_i\partial_iu \cdot\partial_k\varphi_kdx =0 .
\end{aligned}
\end{equation}
 where $\sigma \doteq DH(\varepsilon(u))= \frac{h'(|\varepsilon(u)|)}{|\varepsilon(u)|}\varepsilon(u)$.

Using integration by parts again, we obtain
\ben\label{eq:energy-D2h}
\int_{\Omega} \partial_k\sigma : \varepsilon(\partial_ku) \eta^{2}dx &=& \int_{\Omega} \sigma : \partial_k(D\eta^{2}\otimes \partial_ku)dx +  \int_{\Omega} \partial_k\pi \,\div(\varphi_k)dx  \nonumber\\&&+ \int_{\Omega} u_i\partial_iu \cdot\partial_k\varphi_k \doteq I_1+I_2+I_3,
\een
where $\otimes$ is the tensor product of vectors.

For $I_1$, noticing the relation $|D^2u(x)| \leq 2|D\varepsilon(u)(x)|$, by Young's inequality we have
\begin{equation}
\nonumber
\begin{aligned}
I_1 &= {\int_{\Omega} {\sigma} : {\partial}_k (D\eta^{2}\otimes {\partial}_k u) d x}  \notag \\
&\leq C{\left\{ \int_{\Omega} h'(|\varepsilon(u)|)|Du|( |D{\eta}|^2 + \eta |D^2{\eta}|) d x + \int_{\Omega} h'(|\varepsilon(u)|)|D{\eta}||D^2u|{\eta } dx \right\}}  \notag \\
&\leq {\delta}\int_{\Omega} \frac{h'(|\varepsilon(u)|)}{|\varepsilon(u)|} |D{\varepsilon(u)}|^2\eta^{2} d x + C(\delta)\int_{\Omega} h'(|\varepsilon(u)|)|\varepsilon(u)||D{\eta}|^2  d x  \notag \\
& +C\int_{\Omega} h'(|\varepsilon(u)|)|Du|( |D{\eta}|^2 + \eta |D^2{\eta}|) d x\\
&\leq {\delta} \int_{\Omega} D^2H(\varepsilon(u))(\varepsilon({\partial}_ku),(\varepsilon({\partial}_ku))\eta^{2} d x + C(\delta) \int_{\Omega} h(|\varepsilon(u)|) |D{\eta}|^2 d x \notag \\
& +C\int_{\Omega} h'(|\varepsilon(u)|)|Du|( |D{\eta}|^2 + \eta |D^2{\eta}|) d x.
\end{aligned}
\end{equation}
where $\delta > 0$, to be decided, and we used the estimates (\ref{eq:D2H LOWER BOUND}) and (\ref{eq:th'<h}) in the last step.

For $I_3$, we have
\begin{equation}
\nonumber
\begin{aligned}
I_3 &= \int_{\Omega} u_i{\partial}_i u_j {\partial}_k({\partial}_k u_j \eta^{2}) d x = - \int_{\Omega} {\partial}_k (u_i{\partial}_i u_j){\partial}_k u_j \eta^{2}d x \notag \\
&= - \int_{\Omega} {\partial}_k u_i {\partial}_i u_j {\partial}_k u_j \eta^{2} d x - \int_{\Omega} u_i {\frac12} {\partial}_i(|{\partial}_k u_j|^2) \eta^{2} d x \notag \\
&= {\frac12} \int_{\Omega} |Du|^2u \cdot D\eta^{2} dx,
\end{aligned}
\end{equation}
where we use $\partial_ku_i\partial_iu_j\partial_ku_j = 0$ for divergence free vector $u$ in 2D.

Next, we use (\ref{eq:gns}) to replace $D\pi$, we have
\ben\label{eq:I2}
I_2 &=& \int_{\Omega} \partial_k\pi div(\varphi_k)dx = \int_{\Omega}\partial_k\pi \partial_ku \cdot D\eta^{2} dx \nonumber\\&= & -\int_{\Omega}\sigma_{ik}\partial_i(\partial_ku \cdot D\eta^{2})dx - \int_{\Omega} u_i\partial_iu_k\partial_ku \cdot D\eta^{2}dx.
\een
We estimate (\ref{eq:I2}) in the same way as $I_1$ and $I_3$:
\beno
I_2 & \leq& {\delta} \int_{\Omega} D^2H(\varepsilon(u))(\varepsilon({\partial}_ku),(\varepsilon({\partial}_ku))\eta^{2}d x + C(\delta) \int_{\Omega} h(|\varepsilon(u)|) |D{\eta}|^2 d x \notag \\
& &+C\int_{\Omega} h'(|\varepsilon(u)|)|Du|( |D{\eta}|^2 + \eta |D^2{\eta}|) d x +C\int_\Omega |Du|^2|u||D\eta|\eta dx.
\eeno

Finally, we observe that
\beno
\partial_k\sigma : \varepsilon(\partial_ku)=D^2H(\varepsilon(u))(\varepsilon(\partial_ku),\varepsilon(\partial_ku)).
\eeno
Recall \eqref{eq:energy-D2h} and collect the estimates of $I_1,\cdots,I_3$, and by choosing $\delta$ small enough we deduce
\beno
&&\int_{\Omega}D^2H(\varepsilon(u))(\varepsilon(\partial_ku),\varepsilon(\partial_ku))\eta^{2} dx
 \leq   C \int_{\Omega} h(|\varepsilon(u)|) |D{\eta}|^2 d x \notag \\
&&+C\int_{\Omega} h'(|\varepsilon(u)|)|Du|( |D{\eta}|^2 +\eta |D^2{\eta}|) d x +C\int_\Omega |Du|^2|u||D\eta|\eta dx.
\eeno
Note that  \eqref{eq:D2H LOWER BOUND2}, $|\varepsilon(u)|\leq |D u|$ and (\ref{eq:th'<h}) and we get
\ben\label{eq:energy inequality}
&&\int_{\Omega}|\varepsilon(\partial_ku)|^2\eta^{2} dx
 \leq   C \int_{\Omega} h(|\varepsilon(u)|) |D{\eta}|^2 d x \notag \nonumber\\
&&+C\int_{\Omega} h(|Du|)( |D{\eta}|^2 +\eta |D^2{\eta}|) d x +C\int_\Omega |Du|^2|u||D\eta|\eta dx.
\een

{\bf Step I. The bounded estimate.}

In this step, we choose the cut-off function $\eta=\phi.$
Note that (\ref{eq:phi}) (\ref{eq:energy inequality}), and the energy bounded assumption,
then we deduce that
\beno
\int_{T_{\rho} \backslash T_{2}}|\varepsilon(\partial_ku)|^2dx &\leq& C \frac{1}{(\tau - \rho)^2} +C \frac{1}{\tau - \rho} \int_{T_{\tau} \backslash T_{\rho}\cup (T_3\setminus T_2)} |Du|^2 |u| \phi dx\\
&\doteq& I_1'+I_2'.
\eeno

For the term $I_2'$, noting that
 $\frac{\tau}{2}\leq \rho$ and $u, Du$ is bounded in $T_3\setminus T_2$, we have
\beno
I_2'\leq C \frac{1}{\tau - \rho} \int_{T_{\tau} \backslash T_{\frac{\tau}{2}}} |Du|^2 |u| \phi dx+ C \frac{1}{\tau - \rho}\doteq I_3'+C \frac{1}{\tau - \rho}
\eeno
Next, we deal with the first term of the right hand. Let
\beno
\bar{f}(r)=\frac{1}{2\pi}\int_0^{2\pi}f(r,\theta)d\theta,
\eeno
then by Wirtinger's inequality (for example, for $p=2$ see Ch II.5 \cite{Galdi}) we have
\ben\label{eq:Wirtinger}
\int_0^{2\pi}|f-\bar{f}|^3 \, d\theta\leq C\int_0^{2\pi}|\partial_\theta f|^3d\theta.
\een
By H\"older inequality,
\beno
I_3'\leq C \frac{1}{\tau - \rho} \left(\int_{T_{\tau} \backslash T_{\frac{\tau}{2}}} |Du|^3 \phi dx\right)^{\frac23}\left(\int_{T_{\tau} \backslash T_{\frac{\tau}{2}}}  |u-\bar{u} +\bar{u} |^3 \phi dx\right)^{\frac13}
\eeno

Using (\ref{eq:Wirtinger}) and Lemma \ref{lem:GW} we derive
\ben\label{eq:u L3}
\int_{T_{\tau} \backslash T_{\frac{\tau}{2}}}  |u |^3dx &\leq &\left( \int_{\frac{\tau}{2}<r'< \tau} \int_0^{2\pi}|u(r',\theta)- \bar{u} |^3 \, d\theta \,r' dr' \right)\nonumber\\
&& + \int_{T_{\tau} \backslash T_{\frac{\tau}{2}}} \left( \int_0^{2\pi} u (r,\theta) \, d\theta \right)^3\,dx\nonumber\\
&   \leq&  Cr^3\left( \int_{\frac{\tau}{2}<r'< \tau}  \frac{1}{r'^3}\int_0^{2\pi}|\partial_\theta u|^3d\theta \,r'dr'\right)+C {(\ln r)^{\frac32}r^2}\nonumber\\
&\leq & C r^3\int_{T_{\tau} \backslash T_{\frac{\tau}{2}}}|D u|^3dx+C {(\ln r)^{\frac32}r^2}
\een
for $r>2r_0$($r_0$ is a constant in  Lemma \ref{lem:GW}), since (\ref{eq:h''>0}).


Hence, by using $ (\ln r)^{\frac32}\leq Cr^{\frac12}$ there holds
\beno
I_3'&\leq& C \frac{r}{\tau - \rho} \int_{T_{\tau} \backslash T_{\frac{\tau}{2}}} |Du|^3 dx+Cr^{-\frac12}
\eeno
Recall that the following Poincar\'{e}-Sobolev inequality holds(see, for example, Theorem 8.11 and 8.12 \cite{LL})
\ben
\label{eq:poincare-sobolev}
\|w\|_{L^3(B_\tau\setminus B_{\tau/2})}\leq C \|D w\|_{L^2(B_\tau\setminus B_{\tau/2})}^{\frac13}\|w\|_{L^2(B_\tau)}^{\frac23}+C\tau^{-\frac13}\|w\|_{L^2(B_\tau\setminus B_{\tau/2})},
\een
which implies that
\beno
I_3'&\leq & C \frac{r}{\tau - \rho} \left(\int_{T_{\tau} \backslash T_{\frac{\tau}{2}}} |D^2u|^2dx\right)^{\frac12}\left(\int_{T_{\tau} \backslash T_{\frac{\tau}{2}}} |Du|^2dx\right)+  C \frac{1}{\tau - \rho}\left(\int_{T_{\tau} \backslash T_{\frac{\tau}{2}}} |Du|^2dx\right)^{\frac32} \\
&\leq &\frac1{16}\left(\int_{T_{\tau} \backslash T_{\frac{\tau}{2}}} |D^2u|^2 dx\right)+C \frac{r^2}{(\tau - \rho)^2}
\eeno

Collecting the estimates of $I_1',\cdots,I_3'$, we get
\beno
&&\int_{T_{\rho} \backslash T_{2}}|\varepsilon(\partial_ku)|^2dx \leq \frac1{16}\left(\int_{T_{\tau} \backslash T_{\frac{\tau}{2}}}|D^2u|^2 dx\right)+C \frac{1}{(\tau - \rho)}+C \frac{r^2}{(\tau - \rho)^2}
\eeno
Then by applying $|D^2u(x)| \leq 2|D\varepsilon(u)(x)|$ again,  Lemma \ref{lem:Gia} yields
\beno
\int_{T_{\frac{3}{4} r}\setminus T_2}|D^2 u|^2\,dx\leq C .
\eeno
Finally, by taking $r\rightarrow\infty$, we arrive at
\ben\label{eq:D2u bound}
\int_{\Omega\setminus T_2}|D^2 u|^2\,dx\leq C.
\een

{\bf Step II. The decay estimate.}

In this step, we choose the cut-off function $\eta=\psi.$
Note that (\ref{eq:psi}) (\ref{eq:energy inequality}), and the energy bounded assumption,
then we deduce that
\beno
\int_{T_{2r} \backslash T_{r}}|\varepsilon(\partial_ku)|^2dx &\leq& C \frac{1}{r^2} +C \frac{1}{r} \int_{T_{3r} \backslash T_{r/2}} |Du|^2 |u| \phi dx\\
&\doteq&  C \frac{1}{r^2} +I_4'.
\eeno
Due to $Du\in L^2(\Omega)$ and (\ref{eq:D2u bound}), it follows from Gagliardo-Nirenberg inequality that
\ben\label{eq: Du Lp}
\|Du\|_{L^p(\Omega\setminus{T_3})}<\infty,\quad \forall~p>2.
\een
Thus with the help of Lemma \ref{lem:log growth}, we have
\ben\label{eq:u log}
|u(x)|\leq \sqrt{\ln(|x|)}
\een
for a sufficient large constant, still denoted by $r_0$, and $|x|\geq r_0.$ Consequently,
\beno
\int_{T_{2r} \backslash T_{r}}|\varepsilon(\partial_ku)|^2dx\leq C\frac{\sqrt{\ln r}}{r},
\eeno
which implies the required inequality (\ref{eq:nabla2u decay}).

\end{proof}

\section{Proof of Proposition \ref{prop:nabla3 u}}

\setcounter{equation}{0} \setcounter{theorem}{0}

In this section, we introduce
another cut-off functions, which will be used in the next proof.\\
{\bf Case III. The choosing of $\zeta$.}
First, we introduce
a cut-off function $\zeta \in C_0^{\infty}(\Omega)$ with $ 0 \le \zeta \le 1$, satisfying that\\
i) for $r>2r_0,$ $\rho>0$ and $\tau>0$, there holds $\frac{3}{4}r \le \rho < \tau \le r$;\\
ii)
\begin{equation}\label{eq:eta}\zeta(x)=\zeta(|x|)=
\left\{\begin{array}{llll}
1,\quad ~{\rm in}~T_{2\rho} \backslash T_{r+\frac{r-\rho}{2}};\\
0, \quad ~{\rm in} ~T_{r+\frac{r-\tau}{2}},~{\rm and}~ \Omega\backslash T_{2\tau}
\end{array}\right.
\end{equation}
iii) $|D\zeta| \leq \frac{C}{\tau - \rho},~ |D^2\zeta| \leq \frac{C}{(\tau - \rho)^2}$.

\begin{proof}
For any $\zeta \in C_0^{\infty}(\Omega)$ with $0 \le \zeta \le 1$, letting $\varphi_k = \partial_k\Delta u \,\zeta^3$ with $k=1,2$, we multiply (\ref{eq:gns}) with $\partial_k\varphi_k$ and use integration by parts to obtain $$\int_{\Omega} \partial_k \sigma : \varepsilon(\varphi_k) dx - \int_{\Omega} D\pi \cdot \partial_k \varphi_k dx - \int_{\Omega} u_i\partial_iu \cdot \partial_k \varphi_k dx =0.$$ where $\sigma :=DH(\varepsilon(u)) :=\frac{h^{'}(|\varepsilon(u)|)}{|\varepsilon(u)|} \varepsilon(u).$

Using integration by parts again, by $\varphi_k = \partial_k\Delta u \zeta^3$ we get
\ben\label{eq:energy D3u}
\int_{\Omega} \partial_k \sigma : \varepsilon(\partial_k\Delta u)\zeta^3 dx &=& \int_{\Omega} \sigma : \partial_k(D\zeta^3\otimes \partial_k\Delta u) dx \nonumber\\
&&\mathrel{\phantom{=}}+ \int_{\Omega} \partial_k \pi \div(\varphi_k) dx + \int_{\Omega} u_i\partial_iu \cdot \partial_k \varphi_k dx.
\een
Next we deal with every term of (\ref{eq:energy D3u}):
\beno
&&\int_{\Omega} \partial_k\sigma : \varepsilon(\partial_k \Delta u)\zeta^3 dx \\
 &=& \int_{\Omega} \frac{h^{''}(|\varepsilon(u)|)|\varepsilon(u)| - h^{'}(|\varepsilon(u)|)}{|\varepsilon(u)|^2} \frac{\varepsilon(u) : \varepsilon(\partial_k u)}{|\varepsilon(u)|} \varepsilon(u) : \varepsilon(\partial_k \Delta u) \zeta^3 dx \\
&&\mathrel{\phantom{=}} + \int_{\Omega} \frac{h^{'}(|\varepsilon(u)|)}{|\varepsilon(u)|} \varepsilon(\partial_ku) : \varepsilon(\partial_k \Delta u) \zeta^3 dx  \\
& =&  -\int_{\Omega} \frac{h^{''}(|\varepsilon(u)|)|\varepsilon(u)| - h^{'}(|\varepsilon(u)|)}{|\varepsilon(u)|^2} \frac{\varepsilon(u) : \varepsilon(\partial_{kn} u)}{|\varepsilon(u)|} \varepsilon(u) : \varepsilon(\partial_{kn}u) \zeta^3 dx \\
&&\mathrel{\phantom{=}} - \int_{\Omega} \frac{h^{'}(|\varepsilon(u)|)}{|\varepsilon(u)|} \varepsilon(\partial_{kn}u) : \varepsilon(\partial_{kn} u) \zeta^3 dx \\
&&\mathrel{\phantom{=}} - \int_{\Omega} \partial_n \left(\frac{h^{''}(|\varepsilon(u)|)|\varepsilon(u)| - h^{'}(|\varepsilon(u)|)}{|\varepsilon(u)|^3}\right) (\varepsilon(u) : \varepsilon(\partial_ku)) (\varepsilon(u) : \varepsilon(\partial_{kn}u)) \zeta^3 dx \\
&&\mathrel{\phantom{=}} - \int_{\Omega} \frac{h^{''}(|\varepsilon(u)|)|\varepsilon(u)| - h^{'}(|\varepsilon(u)|)}{|\varepsilon(u)|^2} \frac{\varepsilon(\partial_nu) : \varepsilon(\partial_k u)}{|\varepsilon(u)|} \varepsilon(u) : \varepsilon(\partial_{kn}u) \zeta^3 dx
\\
&&\mathrel{\phantom{=}} - \int_{\Omega} \frac{h^{''}(|\varepsilon(u)|)|\varepsilon(u)| - h^{'}(|\varepsilon(u)|)}{|\varepsilon(u)|^2} \frac{\varepsilon(u) : \varepsilon(\partial_k u)}{|\varepsilon(u)|} \varepsilon(\partial_nu) : \varepsilon(\partial_{kn}u) \zeta^3 dx\\
&&\mathrel{\phantom{=}}  - \int_{\Omega} \partial_n \left(\frac{h^{'}(|\varepsilon(u)|)}{|\varepsilon(u)|} \right) \varepsilon(\partial_ku) : \varepsilon(\partial_{kn}u) \zeta^3 dx  \\
&&\mathrel{\phantom{=}} -\int_{\Omega} \frac{h^{''}(|\varepsilon(u)|)|\varepsilon(u)| - h^{'}(|\varepsilon(u)|)}{|\varepsilon(u)|^2} \frac{\varepsilon(u) : \varepsilon(\partial_{k} u)}{|\varepsilon(u)|} \varepsilon(u) : \varepsilon(\partial_{kn}u) \partial_n(\zeta^3) dx \\
&&\mathrel{\phantom{=}} -  \int_{\Omega} \frac{h^{'}(|\varepsilon(u)|)}{|\varepsilon(u)|} \varepsilon(\partial_ku) : \varepsilon(\partial_{kn}u) \partial_n(\zeta^3) dx \\
&\doteq& -J_1-\cdots-J_8
\eeno

Note that (\ref{eq:D2H LOWER BOUND}), and we have $J_1\geq 0$. By
(\ref{eq;h'' lower bound}), we have
\beno
J_2\geq \int_{\Omega} |\varepsilon(\partial_{kn}u)|^2\zeta^3dx
\eeno
Using $|D^3u(x)| \leq 2|\varepsilon(D^2u)(x)|$, we have
\beno
J_2\geq \frac14\int_{\Omega} |D^3u|^2\zeta^3dx
\eeno

Then
\ben\label{eq:energy}
&&\frac{1}{4}\int_{\Omega} |D^3u|^2\zeta^3dx\leq \int_{\Omega} \partial_k \sigma : (D\zeta^3 \otimes\partial_k\Delta u) dx \nonumber\\
& &- \int_{\Omega} \partial_k \pi div(\varphi_k) dx - \int_{\Omega} u_i\partial_iu \cdot \partial_k \varphi_k dx - J_3-\cdots- J_8.
\een
\par For $J_3$, using (\ref{eq;h''' up bound}), Young's inequality and $|\varepsilon(u)| \leq C$, we have
\beno
J_3&=& \int_{\Omega} \partial_n  \left(\frac{h^{''}(|\varepsilon(u)|)|\varepsilon(u)| - h^{'}(|\varepsilon(u)|)}{|\varepsilon(u)|^3}\right) (\varepsilon(u) : \varepsilon(\partial_ku)) (\varepsilon(u) : \varepsilon(\partial_{kn}u)) \zeta^3 dx \\
 &\leq& C \int_{\Omega} \left(|\varepsilon(u)|^{\gamma_0} + 1 + |\varepsilon(u)|^{\gamma_1 - 1} + |\varepsilon(u)|^{\gamma_2 - 1}\right)|D^2u|^2|D^3u|\zeta^3 dx\\
 &\leq & C \int_{\Omega} |D^2u|^4\zeta^4 dx + \frac{1}{64}\int_{\Omega} |D^3u|^2\zeta^2 dx.
\eeno
where we used
\ben\label{eq:error}
\left|\frac{h^{''}(|\varepsilon(u)|)|\varepsilon(u)| - h^{'}(|\varepsilon(u)|)}{|\varepsilon(u)|^2}\right|\leq C(|\varepsilon(u)|^{\gamma_1 - 1} + |\varepsilon(u)|^{\gamma_2 - 1})
\een
since (\ref{eq:h' lower bound}).

Similarly, $J_6$ can be estimated immediately as
\beno
J_6 &=&\int_{\Omega} \frac{h^{''}(|\varepsilon(u)|)|\varepsilon(u)| - h^{'}(|\varepsilon(u)|)}{|\varepsilon(u)|^2}\partial_n|\varepsilon(u)| \varepsilon(\partial_ku) : \varepsilon(\partial_{kn}u) \zeta^3 dx\\
&\leq& C \int_{\Omega} |D^2u|^4\zeta^4 dx + \frac{1}{64}\int_{\Omega} |D^3u|^2\zeta^2 dx.
\eeno
For $J_4$ and $J_5$, in the same way we have
\beno
J_4 + J_5
& \leq& C\int_{\Omega}\left(|\varepsilon(u)|^{\gamma_2 - 1} + |\varepsilon(u)|^{\gamma_1 - 1}\right) |D^2u|^2 |D^3u| \zeta^3 dx \\
& \leq& C\int_{\Omega}|D^2u|^4 \zeta^4 dx  + \frac{1}{64} \int_{\Omega}|D^3u|^2 \zeta^2 dx.
\eeno

For $J_7$ and $J_8$, using (\ref{eq:error}) and Young inequality again, we get
\beno
J_7+J_8
&\leq & C\int_{\Omega}\left(1+|\varepsilon(u)|^{\gamma_2} + |\varepsilon(u)|^{\gamma_1}\right)|D^2u||D^3u|\zeta^2|D\zeta|dx\\
&\leq & C \int_{\Omega} |D^2u|^2 \zeta^2|D\zeta|^2 dx + \frac{1}{64}\int_{\Omega} |D^3u|^2 \zeta^2dx.
\eeno

Next, we estimate the first three terms of the right hand in (\ref{eq:energy}). Firstly, using (\ref{eq;h'' lower bound}) and (\ref{eq:h' lower bound}) we get
\ben\label{eq1}
A_1&\doteq&\int_{\Omega} \partial_k \sigma : (D\zeta^3 \otimes \partial_k\Delta u) dx\nonumber\\
&= & \int_{\Omega} \partial_k \left(\frac{h'(|\varepsilon(u)|)}{|\varepsilon(u)|}\varepsilon(u)\right): (D\zeta^3 \otimes \partial_k\Delta u) dx \nonumber\\
& \leq &C\int_{\Omega} \left(1 + |\varepsilon(u)|^{\gamma_2}\right) |D^2u|\zeta^2 |D\zeta| |D^3u|dx \nonumber\\
& \leq & C\int_{\Omega} |D^2u|^2 \zeta^2 |D\zeta|^2 dx + \frac{1}{64}\int_{\Omega} |D^3u|^2 \zeta^2dx.
\een
Besides,
\ben\label{eq2}
A_2&\doteq&\int_{\Omega}u_i \partial_iu \cdot \partial_k\varphi_k dx\nonumber\\
&=&-\int_{\Omega} \partial_ku_i \partial_iu_j \partial_k \Delta u_j \zeta^3 dx - \int_{\Omega} u_i \partial_{ki}u_j \partial_k \Delta u_j \zeta^3 dx \nonumber\\
&\leq& C\int_{\Omega} |Du|^2 |D^3u| \zeta^3 dx + C \int_{\Omega} |u| |D^2u| |D^3u| \zeta^3 dx \\
& \leq& \frac{1}{64}\int_{\Omega} |D^3u|^2 \zeta^2 dx + C\int_{\Omega}|u|^2 |D^2u|^2 \zeta^4 dx + C\int_{\Omega} |Du|^4 \zeta^4 dx.\nonumber
\een
Finally, we estimate $A_3=\int_{\Omega} \partial_k \pi \div(\varphi_k) dx$. According the equation of (\ref{eq:gns}),
\begin{equation}
\nonumber
\begin{aligned}
\int_{\Omega} \partial_k \pi \div(\varphi_k) dx &= \int_{\Omega} \partial_k \pi \partial_k \Delta u \cdot D\zeta^3 dx \\&= \int_{\Omega} \partial_i \sigma_{ik} \partial_k \Delta u \cdot D\zeta^3 dx - \int_{\Omega} u_i \partial_iu_k \partial_k \Delta u \cdot D\zeta^3 dx,
\end{aligned}
\end{equation}
whose estimates are similar to $A_1$ and $A_2$, i.e.
\ben\label{eq3}
A_3 & \leq & \frac{1}{64}\int_{\Omega} |D^3u|^2 \zeta^2 dx + C\int_{\Omega} |D^2u|^2 \zeta^2 |D\zeta|^2 dx\nonumber\\
&&+ C\int_{\Omega}|u|^2 |Du|^2|D\zeta|^2 \zeta^2 dx.
\een

Recalling (\ref{eq:energy}), and combining $J_3, \cdots, J_8$, (\ref{eq1}), (\ref{eq2}) and (\ref{eq3}), we have
\beno
\frac{1}{4}\int_{\Omega} |D^3u|^2\zeta^3dx
& \leq & \frac18 \int_{\Omega} |D^3u|^2 \zeta^2 dx + C\int_{\Omega}|u|^2 |D^2u|^2 \zeta^4 dx \nonumber\\
&+& C \int_{\Omega} |Du|^4 \zeta^4 dx + C\int_{\Omega} |D^2u|^2 \zeta^2 |D\zeta|^2dx   \\
&+&  C \int_{\Omega}|D^2u|^4 \zeta^4 dx+ C\int_{\Omega}|u|^2 |Du|^2|D\zeta|^2 \zeta^2 dx .\nonumber
\eeno
Note that the bounded energy integral and (\ref{eq: Du Lp}), i.e. $Du\in L^p(\Omega)$ for any $p\geq 2$, and we derive
\ben\label{eq4}
\frac{1}{4}\int_{\Omega} |D^3u|^2\zeta^3dx
& \leq & \frac18 \int_{\Omega} |D^3u|^2 \zeta^2 dx + C\frac{(\ln r)^{\frac{3}{2}}}{r} \nonumber\\
&+& C+ C\frac{\sqrt{\ln r}}{r(\tau-\rho)^2}+ C\frac{\ln r}{(\tau-\rho)^2}  \nonumber\\
&+&  C \int_{\Omega}|D^2u|^4 \zeta^4 dx
\een
where we also used the decay estimate in Proposition \ref{prop:nabla2 u}. Next we deal with the last term in the right hand.
Due to  Gagliardo-Nirenberg inequality,
\beno
C\int_{\Omega} |D^2u\zeta|^4dx \leq C_0\|D^2u\zeta\|_{L^2(\Omega)}^2 \|D(D^2u\zeta)\|_{L^{2}(\Omega)}^2,
\eeno
where $C_0$ is an absolute constant. Using Proposition \ref{prop:nabla2 u} again, we have
\ben\label{eq5}
\int_{\Omega} |D^2u\zeta|^4dx &\leq& C_0C \frac{\sqrt{\ln r}}{r} \int_{\Omega} |D^3u|^2 \zeta^2 dx+C\frac{\ln r}{r^2(\tau-\rho)^2}\nonumber\\
&\leq &\frac{1}{16}\int_{\Omega} |D^3u|^2 \zeta^2 dx+C\frac{\ln r}{r^2(\tau-\rho)^2}
\een
provided that $r_0$ is large enough.

%

Combining (\ref{eq4}) and (\ref{eq5}), by (\ref{eq:eta}) we deduce
\beno
\int_{T_{2\rho} \setminus T_{r+\frac{r-\rho}{2}}} |D^3u|^2  dx
& \leq & \frac{3}{4} \int_{T_{2\tau} \setminus T_{r+\frac{r-\tau}{2}}} |D^3u|^2  dx + C\frac{(\ln r)^{\frac{3}{2}}}{r} \nonumber\\
&+& C+ C\frac{\sqrt{\ln r}}{r(\tau-\rho)^2}+ C\frac{\ln r}{(\tau-\rho)^2} +C\frac{\ln r}{r^2(\tau-\rho)^2} \nonumber\\
\eeno
Applying Lemma \ref{lem:Gia}, we have
\beno
\int_{T_{2\rho} \setminus T_{r+\frac{r-\rho}{2}}} |D^3u|^2  dx
& \leq & C\frac{(\ln r)^{\frac{3}{2}}}{r} \nonumber\\
&&+ C+ C\frac{\sqrt{\ln r}}{r(\tau-\rho)^2}+ C\frac{\ln r}{(\tau-\rho)^2} +C\frac{\ln r}{r^2(\tau-\rho)^2},
\eeno
by taking $\rho = \frac{3}{4} r$ and $\tau = r$, which implies that
\beno
\int_{T_{\frac{3}{2}r} \setminus T_{\frac{9}{8}r}} |D^3u|^2  dx
&\leq & C\frac{(\ln r)^{\frac{3}{2}}}{r} + C\frac{\sqrt{\ln r}}{r^3} +  \frac{\ln r}{r^2} + C\\
&\leq & C
\eeno
which implies the required inequality (\ref{prop:nabla3 u}).

\end{proof}

\section{Appendix}

\setcounter{equation}{0} \setcounter{theorem}{0}

In the proof of Theorem, we need the following known lemmas.



First, let us recall a result of Gilbarg-Weinberger in \cite{GW1978} about the decay of functions with finite Dirichlet integrals.
\begin{lemma}[Lemma 2.1, 2.2, \cite{GW1978}]
\label{lem:GW}
Let a $C^1$ vector-valued function $f(x)=(f_1,f_2)(x)=f(r,\theta)$ with $r=|x|$ and $x_1=r\cos\theta$.  There holds finite Dirichlet integral in the range $r>r_0$, that is
\beno
\int_{r>r_0}|D f|^2\,dxdy<\infty.
\eeno
Then, we have
\beno
\lim_{r\rightarrow\infty} \frac{1}{\ln r}\int_0^{2\pi}|f(r,\theta)|^2d\theta=0.
\eeno
\end{lemma}


If, furthermore, we assume $D f \in L^p(\R^2)$ for some $2< p < \infty$,
then the above decay property can be improved to be point-wise uniformly.
More precisely, we have

\begin{lemma}[Point-wise Behavior Theorem, Theorem II.9.1 \cite{Galdi}]
\label{lem:log growth}
Let $\Omega \subset \R^2$ be an exterior domain and let
\[
D f \in L^2 \cap L^p (\Omega),
\]
for some $2< p < \infty$.
Then
\[
\lim_{|x| \to \infty} \frac{|f(x)|}{\sqrt{ \ln (|x|)}} = 0,
\]
uniformly.
\end{lemma}

We also need a Giaquinta's iteration lemma \cite[Lemma 3.1]{G83}.

\begin{lemma}[Lemma 3.1 \cite{G83}]
\label{lem:Gia}
Let $f(r)$ be a non-negative bounded function on $[R_0,R_1] \subset \R_+$. If there are negative constants $A,B,D$ and positive exponents $b<a$ and a parameter $\theta \in (0,1)$ such that for all $R_0 \le \rho < \tau \le R_1$
\[
f(\rho) \le \theta f(\tau) + \frac{A}{(\tau-\rho)^{a}} + \frac{B}{(\tau-\rho)^{b}} + D,
\]
then for all $R_0 \le \rho < \tau \le R_1$
\[
f(\rho) \leq C(a,\theta) \left[\frac{A}{(\tau-\rho)^{a}} + \frac{B}{(\tau-\rho)^{b}} + D\right].
\]
\end{lemma}

At last, we introduce the Brezis-Gallouet inequality (see Lemma 2 in \cite{BG1980}, or Lemma 3.1 in \cite{CPZ}).

\begin{lemma}
\label{lem:BG}
Let $f \in H^2(\Omega)$ where $\Omega$ is a bounded
domain or an exterior domain with compact smooth boundary. Then there exists a constant $C_{\Omega}$ depending only on $\Omega$, such that
\beno
\|f\|_{L^{\infty}(\Omega)} \leq C_{\Omega} \|f\|_{H^{1}(\Omega)} \ln^{\frac12} \left(e + \frac{\|D^2f\|_{L^{2}(\Omega)}}{\|f\|_{H^{1}(\Omega)}}\right),
\eeno
or
\beno
\|f\|_{L^{\infty}(\Omega)} \leq C_{\Omega}(1 + \|f\|_{H^{1}(\Omega)}) \ln^{\frac12} \left(e + \|D^2f\|_{L^{2}(\Omega)}\right).
\eeno
\end{lemma}
Note that the second inequality can be obtained immediately from the first one by arguments whether $\|f\|_{H^{1}(\Omega)}<1,$ and we omitted it.

\section*{Acknowledgments}
{ W. Wang was supported by NSFC under grant 11671067 and
 "the Fundamental Research Funds for the Central Universities".
}
\par

\end{document}